\documentclass[12pt]{article}
\usepackage{amssymb}
\usepackage{amsfonts}
\usepackage{amsmath}
\usepackage{graphicx}
\usepackage{float}
\usepackage{geometry}
\usepackage[doublespacing]{setspace}

\newtheorem{theorem}{\sc{Theorem}}[section]
\newtheorem{corollary}[theorem]{\sc{Corollary}}
\newtheorem{definition}[theorem]{\sc{Definition}}
\newtheorem{lemma}[theorem]{\sc{Lemma}}
\newtheorem{example}[theorem]{\sc{Example}}
\newtheorem{remark}[theorem]{\sc{Remark}}
\newenvironment{proof}[1][Proof]{\noindent\textbf{#1.} }{\ \rule{0.5em}{0.5em}}
\input{tcilatex}

\begin{document}

\title{An Alternative Construction to the Transitive Closure of a Directed
Graph}
\author{Kenneth L. Price \\
University of Wisconsin Oshkosh}
\date{\today}
\maketitle

\begin{abstract}
One must add arrows which are forced by transitivity to form the transitive
closure of a directed graph. We introduce a construction of a transitive
directed graph which is formed by adding vertices instead of arrows and
which preserves the transitive relationships formed by distinct vertices in
the original directed graph. \ Our construction does not apply to all
directed graphs.
\end{abstract}

We start by fixing notation and defining compression maps in Section \ref%
{Compression Section}. The motivation for this investigation comes from
algebra through connections with incidence rings, which are rings of
functions defined on sets with relations. None of the algebra applications
are covered here, but we mention anyway a good reference for incidence rings
is \cite{SOD}. The class of generalized incidence rings was introduced by G.
Abrams in \cite{Ab}. His construction is based on the balanced relations. In
Section \ref{Balanced Section} we give the analogous definition for directed
graphs and define stable directed graphs, which form a class between
balanced and preordered directed graphs.

Section \ref{Lemma Section}\ contains our main result, Theorem \ref{Main
Theorem}, which provides conditions for a reflexive directed graph to be the
compression of a preordered directed graph. The proof of Theorem \ref{Main
Theorem} takes up all of section \ref{Proof Section}. Theorem \ref{Main
Theorem} also has a direct application in algebra (see \cite{Price}).

\section{Compression Maps\label{Compression Section}}

The directed graphs we consider have a finite number of vertices and no
repeated arrows. Loops are allowed (a loop is an arrow from a vertex to
itself).\ The vertex set and the arrow set of a directed graph $D$ are
denoted by $V\left( D\right) $ and $A\left( D\right) $, respectively. An
arrow from vertex $v$ to vertex $w$ is denoted by $vw$. The notation $%
D^{\ast }$ is reserved for the subgraph of $D$ with vertex set $V\left(
D^{\ast }\right) =V\left( D\right) $ and arrow set $A\left( D^{\ast }\right)
=\left\{ xy\in A\left( D\right) :x,y\in V\left( D\right) \text{ and }x\neq
y\right\} $.

We say a directed graph $D$ is \emph{reflexive} if $vv\in A\left( D\right) $
for all $v\in V\left( D\right) $ and \emph{transitive} if $xy,yz\in A\left(
D\right) $ implies $xz\in A\left( D\right) $ for all $x,y,z\in V\left(
D\right) $. If $D$ is reflexive and transitive then we say $D$ is \emph{%
preordered}. A \emph{transitive triple}\textit{\ in }$D$ is an ordered
triple of vertices contained in $\limfunc{Trans}\left( D\right) =\left\{
\left( a,b,c\right) :a,b,c\in V\left( D\right) \text{ and }ab,bc,ac\in
A\left( D\right) \right\} $.

\begin{definition}
\label{Compression/Expansion Definition}Suppose $D_{1}$ and $D_{2}$ are
reflexive directed graphs. A \emph{compression map} is a surjective function 
$\theta :V\left( D_{2}\right) \rightarrow V\left( D_{1}\right) $ which
satisfies 1, 2, and 3 below.

\begin{enumerate}
\item $\theta \left( x\right) \theta \left( y\right) \in A\left(
D_{1}\right) $ for all $x,y\in D_{2}$ such that $xy\in A\left( D_{2}\right) $%
.

\item For all $\left( a_{1},a_{2},a_{3}\right) \in \limfunc{Trans}\left(
D_{1}\right) $ there exists $\left( x_{1},x_{2},x_{3}\right) \in \limfunc{%
Trans}\left( D_{2}\right) $ such that $\theta \left( x_{i}\right) =a_{i}$
for $i=1,2,3$.

\item There is a bijection $\theta ^{\ast }:A\left( D_{2}^{\ast }\right)
\rightarrow A\left( D_{1}^{\ast }\right) $ given by $\theta ^{\ast }\left(
xy\right) =\theta \left( x\right) \theta \left( y\right) $ for all $x,y\in
V\left( D_{2}\right) $ with $xy\in A\left( D_{2}^{\ast }\right) $.
\end{enumerate}
\end{definition}

\begin{figure}[th]
\begin{center}
\includegraphics[width=3.14in,height=0.2in]{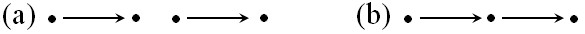}
\end{center}
\caption{(a) is transitive and (b) is not transitive.}
\label{Two Digraphs}
\end{figure}

A figure showing a reflexive directed graph $D$ will only display $D^{\ast }$
and will not show the loops. Thus we assume the directed graphs in Figure %
\ref{Two Digraphs} are both reflexive. In Example \ref{Compress Stable to
Transitive example} we show directed graph (b) is a compression of directed
graph (a).

\begin{example}
\label{Compress Stable to Transitive example}Let $D_{1}$ and $D_{2}$ be the
reflexive directed graphs with $V\left( D_{1}\right) =\left\{ x,y,z\right\} $%
, $A\left( D_{1}^{\ast }\right) =\left\{ xy,yz\right\} $, $V\left(
D_{2}\right) =\left\{ x,y,z,t\right\} $, and $A\left( D_{2}^{\ast }\right)
=\left\{ xy,tz\right\} $ where $x,y,z,t$ are distinct. In Figure \ref{Two
Digraphs} we can match up $D_{1}\ $with (a) and $D_{2}$ with (b). A
compression map $\theta :V\left( D_{2}\right) \rightarrow V\left(
D_{1}\right) $ is given by $\theta \left( x\right) =x$, $\theta \left(
y\right) =y$, $\theta \left( z\right) =z$, and $\theta \left( t\right) =y$.
The effect on the directed graphs is to map the two middle vertices of (a)
to the middle vertex of (b).
\end{example}

\begin{lemma}
\label{transitive inducing lemma}Let $D_{1}$ and $D_{2}$ be reflexive
directed graphs and let $\theta :V\left( D_{2}\right) \rightarrow V\left(
D_{1}\right) \ $be a compression. Suppose $xy,yz\in A\left( D_{2}\right) $
and $\left( \theta \left( x\right) ,\theta \left( y\right) ,\theta \left(
z\right) \right) \in \limfunc{Trans}\left( D_{1}\right) $ for some $x,y,z\in
V\left( D_{2}\right) $. Then $\left( x,y,z\right) \in \limfunc{Trans}\left(
D_{2}\right) $.
\end{lemma}

\proof%
Choose arbitrary $x,y,z\in V\left( D_{2}\right) $ such that $xy,yz\in
A\left( D_{2}\right) $. If $x,y,z$ are not distinct then $xz\in A\left(
D_{2}\right) $ follows immediately. If $x,y,z$ are distinct then set $\theta
\left( x\right) =a$, $\theta \left( y\right) =b$, $\theta \left( z\right) =c$%
. Assume $\left( a,b,c\right) =\left( \theta \left( x\right) ,\theta \left(
y\right) ,\theta \left( z\right) \right) \in \limfunc{Trans}\left(
D_{1}\right) $. Then $a,b,c$ are distinct since $ab,bc,ac\in A\left(
D_{1}^{\ast }\right) $ by part 3 of Definition \ref{Compression/Expansion
Definition}. By part 2 of Definition \ref{Compression/Expansion Definition}
there exist $x^{\prime },y^{\prime },z^{\prime }\in V\left( D_{2}\right) $
such that $\left( x^{\prime },y^{\prime },z^{\prime }\right) \in \limfunc{%
Trans}\left( D\right) $, $\theta \left( x^{\prime }\right) =a$, $\theta
\left( y^{\prime }\right) =b$, and $\theta \left( z^{\prime }\right) =c$.
Moreover $x^{\prime },y^{\prime },z^{\prime }$ are distinct since $a,b,c$
are distinct. Then $x^{\prime }y^{\prime },y^{\prime }z^{\prime }\in A\left(
D_{2}^{\ast }\right) $ so $\theta ^{\ast }\left( xy\right) =ab=\theta ^{\ast
}\left( x^{\prime }y^{\prime }\right) $ and $\theta ^{\ast }\left( yz\right)
=bc=\theta ^{\ast }\left( y^{\prime }z^{\prime }\right) $. Therefore $%
x^{\prime }=x$, $y^{\prime }=y$, and $z^{\prime }=z$ since $\theta ^{\ast }$
is bijective by part 3 of Definition \ref{Compression/Expansion Definition}.
We assumed $\left( x^{\prime },y^{\prime },z^{\prime }\right) \in \limfunc{%
Trans}\left( D\right) \ $and proved $x^{\prime }=x$ and $z^{\prime }=z$ so $%
xz\in A\left( D_{2}\right) .$ Therefore $\left( x,y,z\right) \in \limfunc{%
Trans}\left( D_{2}\right) $ as desired. 
\endproof%

\begin{figure}[th]
\begin{center}
\includegraphics[width=3.62in,height=1.55in]{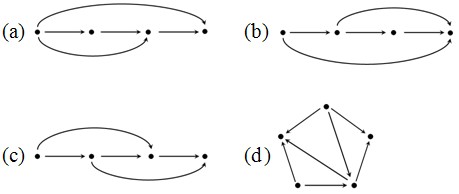}
\end{center}
\caption{Directed graphs which are not compressions of transitive directed
graphs.}
\label{Four Digraphs}
\end{figure}

Lemma \ref{transitive inducing lemma} shows that if a reflexive directed
graph has a preordered compression then it must be a preordered directed
graph.\ Example \ref{Compress Stable to Transitive example} shows a
preordered directed graph may have a compression which is not preordered.
Figure \ref{Four Digraphs} shows there are also directed graphs which are
not compressions of preordered directed graphs. \ Our main result, Theorem %
\ref{Main Theorem}, describes a class of directed graphs which are
compressions of directed graphs.

\section{Balanced and Stable Directed Graphs\label{Balanced Section}}

\begin{definition}
Suppose $D$ is a reflexive directed graph. \ 

\begin{enumerate}
\item $D$ is \emph{balanced} if for all $w,x,y,z\in V\left( D\right) $ such
that $wx,xy,yz,wz\in A\left( D\right) $ there is an arrow from $w$ to $y$ if
and only if there is an arrow from $x$ to $z$.

\item $D$ is \emph{stable} if $D$ is balanced and $ad\in A\left( D\right) $
for all distinct $a,b,c,d\in V\left( D\right) $ such that $ab,ac,bc,bd,cd\in
A\left( D\right) $.
\end{enumerate}
\end{definition}

\begin{figure}[th]
\begin{center}
\includegraphics[width=2.46in,height=0.43in]{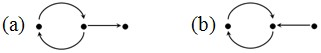}
\end{center}
\caption{Two directed graphs which are not balanced.}
\label{Two Unbalanced}
\end{figure}

A reflexive directed graph is balanced if and only if does not contain an
induced subgraph isomorphic to either (a) or (b) in Figure \ref{Four
Digraphs} or either of the directed graphs in Figure \ref{Two Unbalanced}.

A balanced directed graph is not stable if and only if it does not contain
an induced subgraph isomorphic to (c) in Figure \ref{Four Digraphs}.
Directed graph (d) in Figure \ref{Four Digraphs} is stable but not
preordered.

\begin{theorem}
\label{Compression Theorem}Suppose $\theta :V\left( D_{2}\right) \rightarrow
V\left( D_{1}\right) $ is a compression map and $D_{1}$ and $D_{2}$ are
reflexive directed graphs.

\begin{enumerate}
\item If $D_{1}$ is balanced then $D_{2}$ is balanced.

\item If $D_{1}$ is stable then $D_{2}$ is stable.

\item If $D_{1}$ is preordered then $D_{2}$ is preordered.
\end{enumerate}
\end{theorem}

\proof%
(1) Suppose $D_{1}$ is balanced and $wx,xy,yz,wz\in A\left( D_{2}\right) $
for some $w,x,y,z\in V\left( D_{2}\right) $. Set $a=\theta \left( w\right) $%
, $b=\theta \left( x\right) $, $c=\theta \left( y\right) $, and $d=\theta
\left( z\right) $. Then $ab,bc,cd,ad\in A\left( D_{1}\right) $ by property 1
of Definition \ref{Compression/Expansion Definition}.

If $xz\in A\left( D_{2}\right) $ then property 1 of Definition \ref%
{Compression/Expansion Definition} gives $bd\in A\left( D_{1}\right) $ so $%
ac\in A\left( D_{1}\right) $ since $D_{1}$ is balanced. Then $\left(
a,b,c\right) \in \limfunc{Trans}\left( D_{1}\right) $ so $\left(
w,x,y\right) \in \limfunc{Trans}\left( D_{2}\right) $ by Lemma \ref%
{transitive inducing lemma}. This gives $wy\in A\left( D_{2}\right) $.

For the other direction assume $wy\in A\left( D_{2}\right) $. Then property
1 of Definition \ref{Compression/Expansion Definition} gives $ac\in A\left(
D_{1}\right) $ so $bd\in A\left( D_{1}\right) $ since $D_{1}$ is balanced.
Then $\left( b,c,d\right) \in \limfunc{Trans}\left( D_{1}\right) $ so $%
\left( x,y,z\right) \in \limfunc{Trans}\left( D_{2}\right) $ by Lemma \ref%
{transitive inducing lemma}. This gives $xz\in A\left( D_{2}\right) $.

(2) Suppose $D_{1}$ is stable and $wx,wy,xy,xz,yz\in A\left( D_{2}\right) $
for some distinct vertices $w,x,y,z\in V\left( D_{2}\right) $. Set $a=\theta
\left( w\right) $, $b=\theta \left( x\right) $, $c=\theta \left( y\right) $,
and $d=\theta \left( z\right) $. Then $ab,ac,bc,bd,cd\in A\left( D_{1}^{\ast
}\right) $ by property 3 of Definition \ref{Compression/Expansion Definition}%
. In particular $a,b,c,d$ are distinct and $ad\in A\left( D_{1}^{\ast
}\right) $ since $D_{1}$ is stable. Then $\left( a,c,d\right) \in \limfunc{%
Trans}\left( D_{1}\right) $ so $\left( w,y,z\right) \in \limfunc{Trans}%
\left( D_{2}\right) $ by Lemma \ref{transitive inducing lemma}. This proves $%
wz\in A\left( D_{2}^{\ast }\right) $ as desired. Therefore $D_{2}$ is stable
if $D_{1}$ is stable.

(3) Part 3 follows immediately from Lemma \ref{transitive inducing lemma}.%
\endproof%

The converse does not hold for every part of Theorem \ref{Compression
Theorem}. Figure \ref{Two Unbalanced} shows two directed graphs which are
not balanced. However, they are both compressions of preordered directed
graphs. The compressions are defined in a similar fashion as Example \ref%
{Compress Stable to Transitive example} by constructing a directed graph
which splits the middle vertex in two.

\section{Clasps and Soloists\label{Lemma Section}}

The next definition helps us identify vertices where the transitive relation
fails.

\begin{definition}
Suppose $D$ is a reflexive directed graph and $x\in V\left( D\right) $.

\begin{enumerate}
\item We say $x$ is a \emph{clasp} if there exist $w,y\in V\left( D\right)
\backslash \left\{ x\right\} $ such that $wx,xy\in A\left( D\right) $ and
there is no arrow from $w$ to $y$.

\item We say $x$ is a \emph{locked clasp} if there exist $u,v,w,y\in V\left(
D\right) \backslash \left\{ x\right\} $ such that $\left( u,x,y\right)
,\left( u,x,v\right) ,\left( w,x,v\right) \in \limfunc{Trans}\left( D\right) 
$ and there is no arrow from $w$ to $y$.

\item An \emph{unlocked clasp} is a clasp which is not locked.
\end{enumerate}
\end{definition}

Directed graph (d) in Figure \ref{Four Digraphs} contains a locked clasp
determined by the vertex in the lower right corner. We are now able to state
our main Theorem.

\begin{theorem}
\label{Main Theorem}Let $D$ be a stable directed graph. Then $D$ is the
compression of a preordered directed graph if and only if every clasp in $D$
is unlocked.
\end{theorem}

\begin{remark}
Suppose $D$ is a reflexive directed graph. Theorem \ref{Main Theorem} shows $%
D$ is the compression of a preordered directed graph if $D$ does not contain
an induced subgraph isomorphic to one of the directed graphs in Figure \ref%
{Four Digraphs} or in Figure \ref{Two Unbalanced}. This is reminiscent of
Kuratowki's characterization of planar graphs (see \cite{Kur}). Theorem \ref%
{Main Theorem} is not a complete classification since both directed graphs
in Figure \ref{Two Unbalanced} are compressions of preordered directed
graphs. But both directed graphs in Figure \ref{Two Unbalanced} contain
directed cycles so we can give a complete classification if $D^{\ast }$ is
acyclic.
\end{remark}

\begin{corollary}
Suppose $D$ is a reflexive directed graph such that $D^{\ast }$ is acyclic.
Then $D$ is the compression of a preordered directed graph if and only if $D$
does not contain an induced subgraph isomorphic to one of the directed
graphs in Figure \ref{Four Digraphs}.
\end{corollary}

\begin{definition}
Suppose $D$ is a reflexive directed graph. If $r,s\in V\left( D\right) $
satisfy $rs,sr\in A\left( D\right) $ then $r$ and $s$ are said to be \emph{%
paired} in $V\left( D\right) $. We say an element $s\in V\left( D\right) $
is a \emph{soloist} in $D$ if $r$ and $s$ are not paired for all $r\in
V\left( D\right) \backslash \left\{ s\right\} $.
\end{definition}

\begin{lemma}
\label{Clasp Lemma}Suppose $D$ is a reflexive directed graph.

\begin{enumerate}
\item If $x\in V\left( D\right) $ is a clasp then $x$ is a soloist.

\item Suppose $D_{2}$ is a stable directed graph and $\theta :V\left(
D_{2}\right) \rightarrow V\left( D\right) $ is a compression map. There is a
locked clasp in $D_{2}$ if and only if there is a locked clasp in $D$.

\item If $D$ contains a locked clasp then $D$ is not the compression of a
preordered directed graph.
\end{enumerate}
\end{lemma}

\proof%
(1) Since $x$ is a clasp there exist $w,y\in V\left( D\right) $ such that $%
w,y\in V\left( D\right) \backslash \left\{ x\right\} $, $wx,xy\in A\left(
D\right) $, and there is no arrow from $w$ to $y$. Suppose there exists $%
z\in V\left( D\right) \backslash \left\{ x\right\} $ such that $x$ is paired
with $z$. Applying the balance property to $w,x,z,x$ and to $x,z,x,y$ yields 
$wz,zy\in A\left( D\right) $. This gives $w\neq y$, $y\neq z$, and $w\neq z$
since there is no arrow from $w$ to $y$. Applying the stable property to $%
w,x,z,y$ gives $wy\in A\left( D\right) $, which is a contradiction.

(2) If $x\in V\left( D_{2}\right) $ is a locked clasp then there exist $%
u,v,w,y\in V\left( D_{2}\right) \backslash \left\{ x\right\} $ such that $%
\left( u,x,y\right) ,\left( u,x,v\right) ,\left( w,x,v\right) \in \limfunc{%
Trans}\left( D_{2}\right) $ and there is no arrow from $w$ to $y$ in $D_{2}$%
. Properties (1) and (3) of Definition \ref{Compression/Expansion Definition}
give $\theta \left( u\right) ,\theta \left( v\right) ,\theta \left( w\right)
,\theta \left( y\right) \in V\left( D\right) \backslash \left\{ \theta
\left( x\right) \right\} $ such that $\left( \theta \left( u\right) ,\theta
\left( x\right) ,\theta \left( y\right) \right) ,\left( \theta \left(
u\right) ,\theta \left( x\right) ,\theta \left( v\right) \right) ,\left(
\theta \left( w\right) ,\theta \left( x\right) ,\theta \left( v\right)
\right) \in \limfunc{Trans}\left( D\right) $, and there is no arrow from $%
\theta \left( w\right) $ to $\theta \left( y\right) $ in $D$. Therefore $%
\theta \left( x\right) $ is a locked clasp in $V\left( D\right) $.

If $x_{1}\in V\left( D\right) $ is a locked clasp then there exist $%
u_{1},v_{1},w_{1},y_{1}\in V\left( D\right) \backslash \left\{ x_{1}\right\} 
$ such that $\left( u_{1},x_{1},y_{1}\right) ,\left(
u_{1},x_{1},v_{1}\right) ,\left( w_{1},x_{1},v_{1}\right) \in \limfunc{Trans}%
\left( D\right) $ and there is no arrow from $w_{1}$ to $y_{1}$ in $D$. By
part 2 of Definition \ref{Compression/Expansion Definition} there exist $%
a,b,c,d,e,f,u_{2},x_{2},y_{2}\in V\left( D_{2}\right) $ such that $\left(
a,b,c\right) ,\left( d,e,f\right) ,\left( u_{2},x_{2},y_{2}\right) \in 
\limfunc{Trans}\left( D_{2}\right) $, $\theta \left( a\right) =u_{1}$, $%
\theta \left( b\right) =x_{1}$, $\theta \left( c\right) =v_{1}$, $\theta
\left( d\right) =w_{1}$, $\theta \left( e\right) =x_{1}$, $\theta \left(
f\right) =v_{1}$, $\theta \left( u_{2}\right) =u_{1}$, $\theta \left(
x_{2}\right) =x_{1}$, and $\theta \left( y_{2}\right) =y_{1}$. We have $%
\theta ^{\ast }\left( u_{2}x_{2}\right) =\theta ^{\ast }\left( ab\right) $
and $\theta ^{\ast }\left( bc\right) =\theta ^{\ast }\left( ef\right) $ so $%
u_{2}=a$, $x_{2}=b$, $e=b$, and $c=f$ by part 3 of Definition \ref%
{Compression/Expansion Definition}. Setting $w_{2}=a$ and $v_{2}=c$ gives $%
u_{2},v_{2},w_{2},y_{2}\in V\left( D\right) \backslash \left\{ x_{2}\right\} 
$ such that $\left( u_{2},x_{2},y_{2}\right) ,\left(
u_{2},x_{2},v_{2}\right) ,\left( w_{2},x_{2},v_{2}\right) \in \limfunc{Trans}%
\left( D\right) $, $u_{2}x_{2}$,$x_{2}y_{2}\in A\left( D_{2}\right) $, and
there is no arrow from $u_{2}$ to $y_{2}$ in $D_{2}$. Therefore $x_{2}$ is a
locked clasp in $D_{2}$.

(3) Part 3 follows immediately from part 2.%
\endproof%

\begin{lemma}
\label{Soloist Lemma}Suppose $D$ is a stable directed graph and $s\in
V\left( D\right) $ is a soloist.

\begin{enumerate}
\item Suppose $\left( a,b,s\right) \in \limfunc{Trans}\left( D\right) $ for
some $a,b\in V\left( D\right) \backslash \left\{ s\right\} $ with $a\neq b$.

\begin{enumerate}
\item If $sc\in A\left( D\right) $ for some vertex $c$ then $ac\in A\left(
D\right) $ if and only if $bc\in A\left( D\right) $.

\item If $xa\in A\left( D\right) $ for some vertex $x$ then $xs\in A\left(
D\right) $ if and only if $xb\in A\left( D\right) $.
\end{enumerate}

\item Suppose $\left( a,s,c\right) \in \limfunc{Trans}\left( D\right) $ for
some $a,c\in V\left( D\right) \backslash \left\{ s\right\} $.

\begin{enumerate}
\item If $xa\in A\left( D\right) $ for some vertex $x$ then $xc\in A\left(
D\right) $ if and only if $xs\in A\left( D\right) $.

\item If $cd\in A\left( D\right) $ for some vertex $d$ then $ad\in A\left(
D\right) $ if and only if $sd\in A\left( D\right) $.
\end{enumerate}

\item Suppose $\left( s,b,c\right) \in \limfunc{Trans}\left( V\left(
D\right) \right) $ for some $b,c\in V\left( D\right) \backslash \left\{
s\right\} $ with $b\neq c$.

\begin{enumerate}
\item If $cd\in A\left( D\right) $ for some vertex $d$ then $bd\in A\left(
D\right) $ if and only if $sd\in A\left( D\right) $.

\item If $as\in A\left( D\right) $ for some vertex $a$ then $ab\in A\left(
D\right) $ if and only if $ac\in A\left( D\right) $.
\end{enumerate}
\end{enumerate}
\end{lemma}

\proof%
The proofs of parts 1, 2, and 3 are similar. Note that in part 2 we have $%
a\neq c$ since $s$ is a soloist. We prove part 1. Assume $\left(
a,b,s\right) \in \limfunc{Trans}\left( D\right) $ for some $a,b\in V\left(
D\right) \backslash \left\{ s\right\} $ with $a\neq b$.

(a) Suppose $sc\in A\left( D\right) $ for some $c\in V\left( D\right)
\backslash \left\{ s\right\} $. If $ac\in A\left( D\right) $ then applying
the balance property to $a,b,s,c$ gives $bc\in A\left( D\right) $. On the
other hand if $bc\in A\left( D\right) $ then $a\neq c$ and $b\neq c$ since $%
s $ is a soloist. Applying the stable property to $a,b,s,c$ gives $ac\in
A\left( D\right) $.

(b) Suppose $xa\in A\left( D\right) $ for some $x\in V\left( D\right)
\backslash \left\{ a,b\right\} $. If $xs\in A\left( D\right) $ then applying
the balance property to $x,a,b,s$ gives $xb\in A\left( D\right) $. On the
other hand if $xb\in A\left( D\right) $ then $x\neq s$ since $s$ is a
soloist. Applying the stable property to $x,a,b,s$ gives $xs\in A\left(
D\right) $. 
\endproof%

The proof of Theorem \ref{Main Theorem} is a constructive algorithm
described in section \ref{Proof Section}. In each iteration of the algorithm
we construct a preordered directed graph with one more vertex and define a
compression. The algorithm stops when we arrive at a preordered directed
graph and the desired compression is obtained by composition.

We finish this section with an example which covers the steps and
constructions given in the proof of Theorem \ref{Main Theorem}. The directed
graphs in Figure \ref{Constructions A and B graphs} are stable. We may
identify (i) as a compression of (ii) by mapping 2 and $t_{1}$ to 2. We may
also identify (ii) as a compression of (iii) by mapping 4 and $t_{2}$ to 4.

\begin{example}
Let $D$ be reflexive directed graph (i) shown in Figure \ref{Constructions A
and B graphs}. The clasps are 2 and 4 and we set $x_{1}=2$.

\begin{figure}[th]
\begin{center}
\includegraphics[width=2.1335in,height=1.9813in]{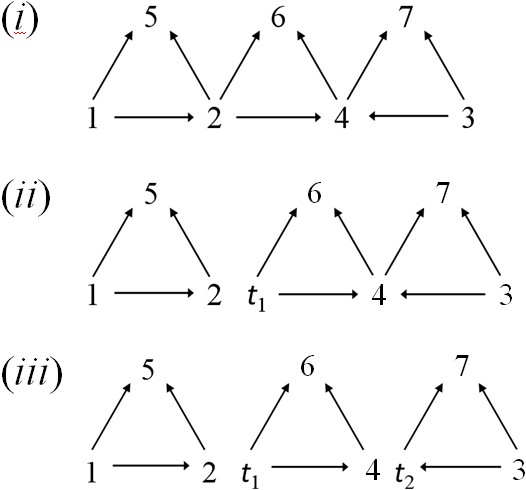}
\end{center}
\caption{A construction using the proof of Theorem \protect\ref{Main Theorem}%
.}
\label{Constructions A and B graphs}
\end{figure}

\begin{description}
\item[Step 1] $Y_{1}=\left\{ 4,6\right\} $ and $A_{1}$ is empty.

\item[Step 2] Use construction A since $A_{1}$ is empty. Let $D_{2}$ be the
reflexive directed graph with $V\left( D_{2}\right) =V\left( D_{1}\right)
\cup \left\{ t_{1}\right\} $ and $A\left( D_{2}^{\ast }\right) =\sigma
_{1}\cup \tau _{1}$ where $B_{1}=\left\{ 4,6\right\} $, $\sigma _{1}=A\left(
D_{1}^{\ast }\right) \backslash \left\{ 24,26\right\} $, and $\tau
_{1}=\left\{ t_{1}4,t_{1}6\right\} $. This gives (ii) in Figure \ref%
{Constructions A and B graphs}.

\item[Step 3] Define $\theta _{1}:V\left( D_{2}\right) \rightarrow V\left(
D_{1}\right) $ by $\theta _{1}\left( t_{1}\right) =2$ and $\theta _{1}\left(
u\right) =u$ for all $u\in V\left( D_{1}\right) $.

\item[Step 4] We go back to step 1 with $x_{2}=4$ since $4$ is the only
clasp in $D_{2}$.

\item[Step 1] We have $Y_{2}=\left\{ 6,7\right\} $ and $A_{2}=\left\{
t_{1},3\right\} $.

\item[Step 2] Use construction B with $a_{2}=3$, $b_{2}=t_{1}$, and $y_{2}=6$%
. Let $D_{3}$ be the reflexive directed graph with $V\left( D_{3}\right)
=V\left( D_{2}\right) \cup \left\{ t_{2}\right\} $ and $A\left( D_{3}^{\ast
}\right) =\sigma _{2}\cup \tau _{2}$ where $\sigma _{2}=A\left( D_{2}^{\ast
}\right) \backslash \left\{ 34,47\right\} $, $\tau _{2}=\left\{
3t_{2},t_{2}7\right\} $. This gives (iii) in Figure \ref{Constructions A and
B graphs}.

\item[Step 3] Define $\theta _{2}:V\left( D_{3}\right) \rightarrow V\left(
D_{2}\right) $ by $\theta _{2}\left( t_{2}\right) =4$ and $\theta _{2}\left(
u\right) =u$ for all $u\in V\left( D_{2}\right) $.

\item[Step 4] $\left( X_{3},\rho _{3}\right) $ is preordered so the
algorithm stops. The compression map is $\theta _{1}\circ \theta _{2}$.
\end{description}
\end{example}

\section{Proof of Theorem \protect\ref{Main Theorem}\label{Proof Section}}

If $D$ is the compression of a preordered directed graph then every clasp in 
$D$ is unlocked by part 3 of Lemma \ref{Clasp Lemma}. We assume $D$ is
stable and every clasp in $D$ is unlocked and prove $D$ is the compression
of a preorder. We set $D_{1}=D$, and describe an algorithm to construct
stable directed graphs $D_{1},\ldots ,D_{m}$ such that $D_{i}$ is a
compression of $D_{i+1}$ for each $i<m$. In the last iteration $D_{m}$ is
preordered and the desired compression map is obtained by composition.

Assume $D_{1}$ is not preordered and fix a clasp $x_{1}\in V\left(
D_{1}\right) $. In the first iteration of our algorithm we have $i=1$.

\subsection{Step 1}

The sets $Y_{i}$ and $A_{i}$ are defined below.

\begin{itemize}
\item $Y_{i}=\left\{ y\in V\left( D_{i}\right) \backslash \left\{
x_{i}\right\} :wx_{i}\text{,}x_{i}y\in A\left( D\right) \text{ and }wy\notin
A\left( D_{i}\right) \text{ for some }w\in V\left( D_{i}\right) \right\} $

\item $A_{i}=\left\{ a\in V\left( D_{i}\right) \backslash \left\{
x_{i}\right\} :\left( a,x_{i},y\right) \in \limfunc{Trans}\left(
D_{i}\right) \text{ for some }y\in Y_{i}\right\} $.
\end{itemize}

Note that $Y_{i}$ is nonempty since $x_{i}$ is a clasp.

\subsection{Step 2}

We fix $t_{i}\notin V\left( D_{i}\right) $ and construct a reflexive
directed graph $D_{i+1}$ such that $V\left( D_{i+1}\right) =V\left(
D_{i}\right) \cup \left\{ t_{i}\right\} $, and $A\left( D_{i+1}^{\ast
}\right) =\sigma _{i}\cup \tau _{i}$ where $\sigma _{i}$ and $\tau _{i}$ are
defined using construction A or construction B. In both constructions $\tau
_{i}=A\left( D_{i+1}^{\ast }\right) \backslash \sigma _{i}$ and $\left\vert
\tau _{i}\right\vert =\left\vert A\left( D_{i}^{\ast }\right) \backslash
\sigma _{i}\right\vert $ so $\left\vert A\left( D_{i}^{\ast }\right)
\right\vert =\left\vert A\left( D_{i+1}^{\ast }\right) \right\vert $. The
arrows in $\tau _{i}$ will all contain $t_{i}$. If an arrow does not contain 
$x_{i}$ then it will be in $\sigma _{i}$. Moreover $\sigma _{i}$ consists of
the arrows belonging to both $D_{i+1}$ and $D_{i}$. Depending on the
construction, an arrow may be contained in $\sigma _{i}$ even if it contains 
$x_{i}$.

Use construction B if there exist $a_{i},b_{i}\in A_{i}$ such that $\left(
b_{i},x_{i},y_{i}\right) \in \limfunc{Trans}\left( D_{i}\right) $ and $%
a_{i}y_{i}\notin A\left( D_{i}\right) $ for some $y_{i}\in Y_{i}$. Otherwise
use construction A.

\subparagraph{Construction A.}

\begin{itemize}
\item $B_{i}=\{b\in V\left( D_{i}\right) \backslash \left\{ x_{i}\right\}
:\left( x_{i},b,y\right) \in \limfunc{Trans}\left( D_{i}\right) $ or $\left(
x_{i},y,b\right) \in \limfunc{Trans}\left( D_{i}\right) $ for some $y\in
Y_{i}\}$

\item $\sigma _{i}=A\left( D_{i}^{\ast }\right) \backslash \left( \left\{
ax_{i}:a\in A_{i}\backslash \left\{ x_{i}\right\} \right\} \cup \left\{
x_{i}b:b\in B_{i}\right\} \right) $

\item $\tau _{i}=\left\{ at_{i}:a\in A_{i}\right\} \cup \left\{ t_{i}b:b\in
B_{i}\right\} $
\end{itemize}

If $y\in Y_{i}$ then $y\in B_{i}$ since $\left( x_{i},y,y\right) \in 
\limfunc{Trans}\left( D_{i}\right) $. Therefore $Y_{i}\subseteq B_{i}$.

\subparagraph{Construction B.}

\begin{itemize}
\item $T_{i}=\left\{ cz:c\in V\left( D_{i}\right) \text{, }z\in V\left(
D_{i}\right) \backslash \left\{ x_{i}\right\} \text{, }\left(
c,x_{i},z\right) \in \limfunc{Trans}\left( V\left( D_{i}\right) \right) 
\text{, and }cy_{i}\notin A\left( D_{i}\right) \right\} $

\item $\sigma _{i}=A\left( D_{i}^{\ast }\right) \backslash \left\{
cx_{i},x_{i}z:c,z\in V\left( D_{i}\right) \text{ and }cz\in T_{i}\right\} $

\item $\tau _{i}=\left\{ ct_{i},t_{i}z:c,z\in V\left( D_{i}\right) \text{
and }cz\in T_{i}\right\} $
\end{itemize}

There exists $z\in Y_{i}$ such that $\left( a_{i},x_{i},z\right) \in 
\limfunc{Trans}\left( D_{i}\right) $ since $a_{i}\in A_{i}$. Moreover $%
a_{i}y_{i}\notin A\left( D_{i}\right) $, $z\in x_{i}$, and $x_{i}\notin
Y_{i} $ so $a_{i}z\in T_{i}$. Therefore $T_{i}$ is nonempty.

\subsection{Step 3}

Define $\theta _{i}:V\left( D_{i+1}\right) \rightarrow V\left( D_{i}\right) $
so that $\theta _{i}\left( t_{i}\right) =x_{i}$ and $\theta _{i}\left(
u\right) =u$ for all $u\in V\left( D_{i}\right) $.

Before moving on we prove $\theta _{i}$ is a compression. Routine
calculations show part (1) of Definition \ref{Compression/Expansion
Definition} hold and there is a well-defined map $\theta _{i}^{\ast
}:A\left( D_{i+1}^{\ast }\right) \rightarrow A\left( D_{i}^{\ast }\right) $
given by $\theta _{i}^{\ast }\left( uv\right) =\theta _{i}\left( u\right)
\theta _{i}\left( v\right) $ for all $u,v\in V\left( D_{i}\right) $ such
that $uv\in A\left( D_{i+1}^{\ast }\right) $. It is easy to see $\theta
_{i}^{\ast }\left( \sigma _{i}\right) \cup \theta _{i}^{\ast }\left( \tau
_{i}\right) =A\left( D_{i}^{\ast }\right) $ hence $\theta _{i}^{\ast }$ is
surjective. We have already shown $\left\vert A\left( D_{i+1}^{\ast }\right)
\right\vert =\left\vert A\left( D_{i}^{\ast }\right) \right\vert $ so $%
\theta _{i}^{\ast }$ is a bijection.

The only condition left is part 2 of Definition \ref{Compression/Expansion
Definition}. Suppose $d_{1},d_{2},d_{3}\in V\left( D_{i}\right) $ and $%
\left( d_{1},d_{2},d_{3}\right) \in \limfunc{Trans}\left( D_{i}\right) $. We
must show $\left( d_{1},d_{2},d_{3}\right) $ is the image of a transitive
triple in $D_{i+1}$. This is easy if $d_{1},d_{2},d_{3}$ are not distinct
since every arrow in $D_{i}$ is the image of an arrow in $D_{i+1}$. Assume $%
d_{1},d_{2},d_{3}$ are distinct.

We check every possible case and make\ repeated use of the fact $uv\in
\sigma _{i}$ if and only if $uv\in A\left( D_{i}^{\ast }\right) $ for all $%
u,v\in V\left( D_{i}\right) \backslash \left\{ x_{i}\right\} $. In cases 2,
3, and 4 we have $x_{i}=d_{r}$ for some $r\in \left\{ 1,2,3\right\} $. We
will show\ either $\left( d_{1},d_{2},d_{3}\right) \in \limfunc{Trans}\left(
D_{i}\right) $ so that $\theta _{i}\left( d_{j}\right) =d_{j}$ for $j=1,2,3$
or the desired transitive triple is obtained by replacing $d_{r}$ with $%
t_{i} $ so that $\theta _{i}\left( d_{j}\right) =d_{j}$ for $j\neq r$, and $%
\theta _{i}\left( t_{i}\right) =d_{r}$.

We split cases between construction B and construction A when necessary.
Note that $x_{i}$ is a soloist by part 1 of Lemma \ref{Clasp Lemma}.

\begin{description}
\item[Case 1] $d_{1},d_{2},d_{3}\in V\left( D_{i}\right) \backslash \left\{
x_{i}\right\} $
\end{description}

We have $d_{1}d_{2},d_{2}d_{3},d_{1}d_{3}\in \sigma _{i}$ thus $\left(
d_{1},d_{2},d_{3}\right) \in \limfunc{Trans}\left( D_{i}\right) $.

\begin{description}
\item[Case 2] If $d_{1},d_{2}\in V\left( D_{i}\right) \backslash \left\{
x_{i}\right\} $ and $d_{3}=x_{i}$ then $d_{1}d_{2}\in \sigma _{i}$ since $%
d_{1}d_{2}\in A\left( D_{i}^{\ast }\right) $.
\end{description}

\subparagraph{Check case 2 for construction A.}

If $x_{i}y\in A\left( D_{i}\right) $ then $d_{1}y\in A\left( D_{i}\right) $
if and only if $d_{2}y\in A\left( D_{i}\right) $ by part 1(a) of Lemma \ref%
{Soloist Lemma}. This gives $d_{1}\in A_{i}$ if and only if $d_{2}\in A_{i}$
so $d_{1}t_{i}\in \tau _{i}$ if and only if $d_{2}t_{i}\in \tau _{i}$ and
either $\left( d_{1},d_{2},t_{i}\right) \in \limfunc{Trans}\left(
D_{i}\right) $ or $\left( d_{1},d_{2},x_{i}\right) \in \limfunc{Trans}\left(
D_{i}\right) $.

\subparagraph{Check case 2 for construction B.}

We have $d_{1}z\in A\left( D_{i}\right) $ if and only if $d_{2}z\in A\left(
D_{i}\right) $ for all $z\in V\left( D_{i}\right) \backslash \left\{
x_{i}\right\} $ such that $x_{i}z\in A\left( D_{i}\right) $ by part 1(a) of
Lemma \ref{Soloist Lemma}. This gives $d_{1}z\in T_{i}$ if and only if $%
d_{2}z\in T_{i}$ for all $z\in V\left( D_{i}\right) \backslash \left\{
x_{i}\right\} $. Therefore $d_{1}t_{i}\in \tau _{i}$ if and only if $%
d_{2}t_{i}\in \tau _{i}$ and either $\left( d_{1},d_{2},t_{i}\right) \in 
\limfunc{Trans}\left( D_{i}\right) $ or $\left( d_{1},d_{2},x_{i}\right) \in 
\limfunc{Trans}\left( D_{i}\right) $.

\begin{description}
\item[Case 3] If $d_{1},d_{3}\in V\left( D_{i}\right) \backslash \left\{
x_{i}\right\} $ and $d_{2}=x_{i}$ then $d_{1}d_{3}\in \sigma _{i}$ since $%
d_{1}d_{3}\in A\left( D_{i}^{\ast }\right) $.
\end{description}

\subparagraph{Check case 3 for construction A.}

If $d_{3}\in B_{i}$ then there exists $y\in Y_{i}$ such that $\left(
x_{i},d_{3},y\right) \in \limfunc{Trans}\left( D_{i}\right) $ or $\left(
x_{i},y,d_{3}\right) \in \limfunc{Trans}\left( D_{i}\right) $. This gives $%
d_{1}y\in A\left( D_{i}\right) $ by applying either part 2(b) or part 3(b)
of Lemma \ref{Soloist Lemma}. Therefore $d_{1}\in A_{i}$.

On the other hand if $d_{1}\in A_{i}$ then $\left( d_{1},x_{i},z\right) \in 
\limfunc{Trans}\left( D_{i}\right) $ for some $z\in Y_{i}$. Since $z\in
Y_{i} $ there exists $w\in V\left( D_{i}\right) $ such that $wx_{i}\in
A\left( D_{i}\right) $ and $wz\notin A\left( D_{i}\right) $. If $wd_{3}\in
A\left( D_{i}\right) $ then $\left( d_{1},x_{i},z\right) ,\left(
d_{1},x_{i},d_{3}\right) ,\left( w,x_{i},d_{3}\right) \in \limfunc{Trans}%
\left( D_{i}\right) $ and $x_{i}$ is a locked clasp, which is a
contradiction. We are left with $wd_{3}\notin A\left( D_{i}\right) $, $%
d_{3}\in Y_{i}$, and $d_{3}\in B_{i}$.

We have shown $d_{1}\in A_{i}$ if and only if $d_{2}\in B_{i}$ so $%
d_{1}t_{i}\in \tau _{i}$ if and only if $t_{i}d_{3}\in \tau _{i}$ and either 
$\left( d_{1},x_{i},d_{3}\right) \in \limfunc{Trans}\left( D_{i}\right) $ or 
$\left( d_{1},t_{i},d_{3}\right) \in \limfunc{Trans}\left( D_{i}\right) $.

\subparagraph{Check case 3 for construction B.}

If $d_{1}y_{i}\notin A\left( D_{i}\right) $ then $d_{1}d_{3}\in T_{i}$ and $%
\left( d_{1},t_{i},d_{3}\right) \in \limfunc{Trans}\left( D_{i}\right) $.

If $d_{1}y_{i}\in A\left( D_{i}\right) $ then $d_{1}x_{i}\in \sigma _{i}$
and we must show $x_{i}d_{3}\in \sigma _{i}$. Note that $\left(
d_{1},x_{i},d_{3}\right) ,\left( d_{1},x_{i},y_{i}\right) ,\left(
b_{i},x_{i},y_{i}\right) \in \limfunc{Trans}\left( D_{i}\right) $ so $%
b_{i}d_{3}\in A\left( D_{i}\right) $ since $x_{i}$ is an unlocked clasp. If $%
x_{i}d_{3}\notin \sigma _{i}$ then $cd_{3}\in T_{i}$ for some $c\in V\left(
D_{i}\right) $ such that $cy_{i}\notin A\left( D_{i}\right) $. This gives $%
\left( b_{i},x_{i},y_{i}\right) ,\left( b_{i},x_{i},d_{3}\right) ,\left(
c,x_{i},d_{3}\right) \in \limfunc{Trans}\left( D_{i}\right) $ with $%
cy_{i}\notin A\left( D_{i}\right) $ and $x_{i}$ is a locked clasp. This is a
contradiction.\ We are left with $d_{1}x_{i},x_{i}d_{3}\in \sigma _{i}$ and $%
\left( d_{1},x_{i},d_{3}\right) \in \limfunc{Trans}\left( D_{i}\right) $.

\begin{description}
\item[Case 4] If $d_{2},d_{3}\in V\left( D_{i}\right) \backslash \left\{
x_{i}\right\} $, $d_{1}=x_{i}$ then $d_{2}d_{3}\in \sigma _{i}$ since $%
d_{2}d_{3}\in A\left( D_{i}^{\ast }\right) $.
\end{description}

\subparagraph{Check case 4 for construction A.}

Suppose $d_{j}\in B_{i}$ for $j=2$ or $j=3$ and let $k\in \left\{
2,3\right\} $ be such that $k\neq j$. If $d_{j}\in B_{i}$ then $\left(
x_{i},d_{j},y\right) \in \limfunc{Trans}\left( D_{i}\right) $ or $\left(
x_{i},y,d_{j}\right) \in \limfunc{Trans}\left( D_{i}\right) $ for some $y\in
Y_{i}$. There exists $w\in V\left( D_{i}\right) \backslash \left\{
x_{i}\right\} $ such that $wx_{i}\in A\left( D_{i}\right) $ and $wy\notin
A\left( D_{i}\right) $ since $y\in Y_{i}$. Then $wd_{2},wd_{3}\notin A\left(
D_{i}\right) $ by two applications of part 3(b) of Lemma \ref{Soloist Lemma}%
. Thus $d_{2},d_{3}\in Y_{i}$ and $\left( t_{i},d_{2},d_{3}\right) \in 
\limfunc{Trans}\left( D_{i}\right) $.

We have shown $d_{2}\in B_{i}$ or $d_{3}\in B_{i}$ imply $\left(
t_{i},d_{2},d_{3}\right) \in \limfunc{Trans}\left( D_{i}\right) $. On the
other hand if $d_{2}\notin B_{i}$ and $d_{3}\notin B_{i}$ then $%
x_{i}d_{2},x_{i}d_{3}\in \sigma _{i}$ and $\left( x_{i},d_{2},d_{3}\right)
\in \limfunc{Trans}\left( D_{i}\right) $.

\subparagraph{Check case 4 for construction B.}

We have $cd_{2}\in A\left( D_{i}\right) $ if and only if $cd_{3}\in A\left(
D_{i}\right) $ for all $c\in V\left( D_{i}\right) $ such that $cx_{i}\in
A\left( D_{i}\right) $ and $cy_{i}\notin A\left( D_{i}\right) $ by part 3(b)
of Lemma \ref{Soloist Lemma}. This gives $cd_{2}\in T_{i}$ if and only if $%
cd_{3}\in T_{i}$ for all $c\in V\left( D_{i}\right) \backslash \left\{
x_{i}\right\} $. Therefore $t_{i}d_{2}\in \tau _{i}$ if and only if $%
t_{i}d_{3}\in \tau _{i}$ and either $\left( x_{i},d_{2},d_{3}\right) \in 
\limfunc{Trans}\left( D_{i}\right) $ or $\left( t_{i},d_{2},d_{3}\right) \in 
\limfunc{Trans}\left( D_{i}\right) $.

\subsection{Step 4}

If $D_{i+1}$ is preordered then the algorithm stops and the compression from 
$D_{i+1}$ to $D$ is determined by composition. Otherwise fix a clasp $%
x_{i+1}\in V\left( D_{i}\right) $ and\ go back to step 1 with $i$ replaced
by $i+1$.

To study the algorithm we consider a given iteration $i$. Then $A\left(
D_{i+1}\right) $ is stable by Theorem \ref{Compression Theorem} and $D_{i+1}$
contains no unlocked clasps\ by part 2 of Lemma \ref{Clasp Lemma}. This
means we may repeat the algorithm as often as necessary. We must prove the
algorithm stops eventually.

In each iteration of the algorithm we are adding a new vertex but not adding
any arrows other than loops. The only way this can continue indefinitely is
if our algorithm forces vertices to not form arrows with any other elements.
We assume every vertex of $D_{i}$ forms an arrow with some other vertex of $%
D_{i}$ and show every vertex of $D_{i+1}$ forms an arrow with some other
vertex of $D_{i+1}$.

Suppose $x,y\in V\left( D_{i}\right) $ satisfy $x\neq y$ and $xy\in A\left(
D_{i}\right) $. If $x\neq x_{i}$ and $y\neq x_{i}$ then $xy\in \sigma _{i}$
so $xy\in A\left( D_{i+1}\right) $. Note that $t_{i}$ forms an arrow with
some other vertex of $D_{i+1}$ by construction. We are left with proving $%
x_{i}$ forms an arrow with some other vertex of $D_{i+1}$.

Assume there is not an arrow formed by $x_{i}$ and any other vertex of $%
D_{i+1}$ after using construction A. There exist $b,z\in V\left(
D_{i}\right) $ such that $bx_{i}\in A\left( D_{i}\right) $, $x_{i}z\in
A\left( D_{i}\right) $, and $bz\notin A\left( D_{i}\right) $ since $x_{i}$
is a clasp. Then $bt_{i}\in A\left( D_{i+1}\right) $ and $t_{i}z\in A\left(
D_{i+1}\right) $ by assumption so $b\in A_{i}$ and $z\in B_{i}$. Since $b\in
A_{i}$ there exists $y\in V\left( D_{i}\right) $ such that $\left(
b,x_{i},y\right) \in \limfunc{Trans}\left( D_{i}\right) $. There must also
exist $a\in V\left( D_{i}\right) $ such that $ax_{i}\in A\left( D_{i}\right) 
$ and $ay\notin A\left( D_{i}\right) $ since $y\in Y_{i}$. This gives $%
at_{i}\in A\left( D_{i+1}\right) $ by assumption so $a\in A_{i}$. Thus $%
a,b\in A_{i}$ satisfy the conditions in step 2 for construction B. This
contradicts our assumption that we used construction A, so there is an arrow
formed by $x_{i}$ and another vertex $\ $of $D_{i+1}$.

In construction B we have $b_{i}y_{i}\in A\left( D_{i}\right) $ so $%
b_{i}z\notin T_{i}$ for all $z\in V\left( D_{i}\right) \backslash \left\{
x_{i}\right\} $. This gives $b_{i}x_{i}\in \sigma _{i}$ and $b_{i}x_{i}\in
A\left( D_{i+1}\right) $ so there is an arrow formed by $x_{i}$ with another
vertex of $D_{i+1}$.

\paragraph{Acknowledgements}

The author wishes to thank Martin Erickson of Truman State University and
Steven Szydlik of the University of Wisconsin Oshkosh for encouragement and
suggesting improvements to the writing in earlier drafts.

\end{document}